\documentclass[12pt]{article}
\usepackage{amsmath}
\usepackage{amssymb}
\usepackage{url}

\usepackage[T1]{fontenc}
\usepackage{charter}
\usepackage{diagrams}
\usepackage{tikz, tikz-cd}

\diagramstyle{midshaft,PostScript=dvips,nohug}

\newtheorem{lemma}{Lemma}[section]
\newtheorem{prop}[lemma]{Proposition}

\newtheorem{defn}[lemma]{Definition}

\newcommand{\E}{{\mathcal E}}
\newcommand{\X}{{\mathcal X}}

\newcommand{\M}{{\mathfrak M}}

\newcommand{\J}{{\mathbb J}}

\newcommand{\imes}{\times}
\newcommand{\gd}{(generalized) }

\newarrow{Line}{-}{-}{-}{-}{-}
\newarrow{Dashto}{}{dash}{}{dash}>
\newarrow{Corresponds}{<}{-}{-}{-}{>}

\title{Elementwise semantics in categories with pull-backs}
\author{Anders Kock}
\date{}

\begin{document}

\maketitle

\section*{Introduction}In several places, including \cite{[SDG]}, II,1-3, is expounded how to work with \gd 
{\em elements} of the objects in a 
category $\E$ with finite limits (``Krip\-ke-Joyal semantics''). 

In the present note, we describe and 
extend this ``semantics'' when the category $\E$  just is assumed to have 
pull-backs, thus no terminal object $1$ is assumed, nor binary products 
$A \times B$.
\footnote{The reason for not insisting in the finite products (in 
particular, a terminal object $1$) is 
the possible application in e.g.\ the partial toposes, as studied in 
\cite{[BS]}. 
Note that there is a qualitative difference between the 
categorical properties $\E$ itself and the properties of its 
``slices'' $\E/X$: for $\E/X$ always has a terminal object; this is 
not assumed for $\E$.}

The basic notions are

\medskip

- \gd element $a$ of an object $A$ of $\E$

- \gd subobject\footnote{Unlike most standard references on 
Kripke-Joyal semantics like \cite{[MM]} VI.6, we do not assume that $\E$ is a 
topos. In a topos, {\em subobjects} of $A$ may be encoded as  
{\em elements} of $\Omega^{A}$; and the semantics for subobjects and 
partial maps may be 
reduced to semantics for elements. The same applies to the \gd 
subobjects and \gd partial maps.}
 $U$ of an object $A$ of $\E$

- \gd partial map $s$ from an object $A$ of $\E$ to another object $E$ of 
$\E$.

Not included here is the elementwise semantics ``object of maps'' from $A$ to $B$, 
(which will be an object $B^{A}$) and related constructs, possible in 
Cartesian closed categories or related kinds of categories. 
Elementwise semantics for certain such constructs exist, see e.g.\ \cite{[SDG]} 
II.4. We do consider such ``higher order'' objects here, but only with traditional 
diagram-style tools, see the section on jet {\em bundles} in 
\ref{Jebux} below. 
We intend in a later note also to give 
such things an elementwise semantics.

 The maps in $\E$ will  be called ``actual maps'', or 
even just ``maps'' or ``arrows''.

I was led to the desire for an explicit semantics for these things, 
because it is in this  way that some synthetic reasoning, e.g.\ for the 
geometry of jet bundles, as in \cite{[SGM]} 22.7 or in \cite{[BFF]} Theorem 11.1, may be fully justified and 
communicated. We deal with a theory of (section-) jets below; but I 
believe that the methodology we develop for it in the present note 
has a more general scope.

\medskip

Some derived relations 
\medskip

- $\in$: when is a \gd element $a$ of $A$ {\em contained in} or a {\em member 
of} a \gd subobject $U$ of  $A$; 

- $\subseteq $: the partial order of  \gd subobjects of $A$

- \gd partial map

- {\em support} of a  \gd partial map from $A$ to $E$,  as a  \gd subobject  of $A$

- {\em value} $s(a)$ of a  \gd partial map on a  \gd  element $a$  which is a member 
of the support of $s$.

- {\em counterimage} of a  \gd subobject of $A$ along a map $A'\to A$.

\medskip

And finally, we have  some results (``principles''), e.g.\ 
the extensionality principles:

\medskip

- a  \gd subobject is determined by the  \gd elements contained in it

- a  \gd partial map $s$ is determined by the values of $s$ on the 
 \gd elements contained in its support.

\medskip

\sloppy The \gd elements, \gd subobjects, \gd partial maps,  are all given at a 
certain ``stage'' $X\in \E$, which is subject to variation (``change of 
state'') as in Kripke's work, and explained in category theoretic 
terms in \cite{[SDG]}, say; a stage $X$ is  any object of 
$\E$, and the ``change of stage'' from $X$ to $Y$ takes place along a 
map $\alpha: Y\to X$ in $\E$. (This change of stage can be formulated by 
saying that the 
notions etc.\ described are contravariant, set-valued constructions, 
i.e.\ that they take place in the functor category $Set ^{\E ^{op}}$. 
More generally, the change of stage can be encoded as passing from the 
slice category $\E/X$ to $\E/Y$ by ``the'' pull-back functor $\alpha 
^{*}: \E/X\to \E/Y$. We avoid here using this useful technique, because it involves the 
choice of definite 
pull-backs, - but see e.g.\ \cite{[SDG]} for an exposition of this slice 
category  technique. The uses that we shall make of the notation $\alpha^{*}$ 
are {\em exact}, i.e.\ do not depend on any choice, or the coherence issues 
that arise from resulting comparison isomorphisms.

If $\E$ happens to have finite (chosen) products, a \gd subobject $U$ of $A$ 
at stage $X$ can be encoded as an actual subobject of $A\times X$.

Not included in the present note is the consideration of generalized {\em 
map} from $A$ to $B$; a generalized map at stage $X$ is an actual map 
$A\imes X \to B$ (assuming binary products in $\E$).

The decoration of an entity, like a \gd subobject $U$ of $A$ at stage 
$X$, will be $U\subseteq _{X}A$; the subscript indicates that we are 
talking about a {\em generalized} subobject defined at stage $X$, not about an 
actual subobject of $U\subseteq A$, in the 
standard sense of the category $\E$. If $\E$ happens to have a 
terminal object $1$, $U\subseteq _{1}A$ will be equivalent to  
$U\subseteq A$ in the standard sense.

In the sequel, the phrase ``generalized'' is sometimes omitted.

\section{Basic notions}

 \subsection{Elements: $a\in_{X}A$}\label{ss1} If $A$ is an object of $\E$, 
and $a:X\to A$ is an arrow in $\E$, we say that $a\in_{X}A$, or that $a$ 
is a {\em generalized} element of $A$ {\em defined at stage} $X$, or {\em parametrized} 
by $X$. 
If $f:A \to B$ is an actual map in $\E$ and $a\in_{X}A$, we have $f\circ a \in 
_{X}B$, and it may also be denoted $f(a)$. This usage will be extended to 
\gd partial maps from $A$ to $B$. Note that we compose maps from the 
right to the left.

If 
$\alpha :Y \to X$, we define $\alpha^{*}(a)\in _{Y}A$ to be $a\circ 
\alpha$, or ``$a$ considered at the later stage $Y$''. The 
associative law for $\circ$ implies that $\alpha^{*}(f(a))= 
f(\alpha^{*}(a))$; both equal $f\circ a\circ \alpha$.

 Consider 
a pull-back square obtained from two given maps $f$ and $p$ with 
common codomain, as in
\begin{equation}\label{pb1x}\begin{diagram}M \SEpbk & \rTo ^{d} &B\\
\dTo^{c} &&\dTo_{p}\\
A&\rTo_{f}&C
\end{diagram}\end{equation}

Then a \gd  element at stage $X$ of $M$ is uniquely 
given by a {\em pair} of \gd elements: $a\in _{X} A$ and $b\in _{X}B$ 
with the property that $f(a)= p(b) \in_{X}A$. This is just a 
reformulation of the universal property of a pull-back; and it does not 
mention any other data than the two given maps $f$ and $p$.
The unique \gd element of $M$ thus given, one denotes 
$\langle a, b \rangle \in _{X}M$.

\subsection{Subobjects and generalized subobjects}\label{Subox}

Recall that a {\em subobject} $U$ of an object 
$A$ in a category $\E$ is {\em represented by} a monic map $i: 
N\rightarrowtail A$; and $i': N'\rightarrowtail A$ represents the same 
subobject if there is some (necessarily unique) isomorphism $\mu : N\to 
N'$ with $i'\circ \mu  =i$, cf.\ \cite{[M]} V.7; thus $U$ is an equivalence 
class of monics with codomain $A$. We want to be pedantic about 
the distinction between a subobject and a representing monic. We 
write $U\subseteq A$ when $U$ is a subobject of $A\in \E$.  But note that a 
subobject $U$ of $A$ is not an object of $\E$. We shall use similar 
pedantry also for some other of the following notions. 

A subobject $U$ of $A$ in 
this sense, we shall also call an {\em actual} subobject of $A$, to 
distinguish it from a {\em generalized} subobject, {\em defined 
at stage} $X\in \E$, as in Definition \ref{AAAx} below.

Given objects $A$ and $B$ in $\E$.
\begin{defn}\label{Relx}A {\em relation} $\M$ from  $A$ to 
 $B$  is an 
 equivalence class  of 
 jointly monic spans
 \begin{equation}\label{relx}\begin{diagram}A&\lTo^{c}&M&\rTo^{d}&B
 \end{diagram}.\end{equation}
 \end{defn}
 
 The equivalence relation is the evident one: $(c,d)$ is 
 equivalent to $(c',d')$ if there is an isomorphism $\mu: M\to M'$ (necessarily 
 unique)  with $c'\circ \mu =c$ and $d'\circ \mu=d$.
 
 If binary Cartesian products are available, and chosen, in $\E$,  
 relations from $A$ to $B$ may be identified with subobjects of $A \times 
 B$.
 
 Given two jointly monic spans, as 
displayed with full 
arrows in
\begin{equation}\label{mux}\begin{diagram}A&\lTo^{c}&M&\rTo^{d}&B\\
&\luTo_{c'}&\dDashto _{\mu}&\ruTo_{d'}\\
&&M'&&
\end{diagram}\end{equation}
then because $c',d'$ are jointly monic, there is at most one $\mu:M 
\to M'$ (dotted arrow) making the two triangles commute. If there is 
such a $\mu$, we say $(c,d)\leq (c',d')$; this defines a preorder on the class 
of jointly monic spans from $A$ to $B$.   
We write $\M \subseteq \M'$ if $(c,d)\leq (c',d')$, where $\M$ is the 
relation 
represented by $(c,d)$ and similarly $\M'$ is represented by $(c',d')$.

 Pull-backs provide examples of  relations: give a pair $f,p$ of arrows 
 with common codomain, as in (\ref{pb1x}), the set of pairs of arrows 
 $c,d$ completing $f,p$ into a pull-back square is a relation, i.e.\ 
 an equivalence class of jointly monic spans.
 
 We now consider a relation $\M$ from $A$ to $B$ in its 
 contravariant dependence of $B$, which we therefore denote $X$:

Given objects $A$ and $X$ in $\E$. With the same equivalence relation 
as in Definition \ref{Relx} (with $X$ for $B$), we pose:
 \begin{defn}\label{AAAx}A (generalized) subobject $U$ of $A$ at 
 stage $X$ (written $U\subseteq_{X}A$) is an 
 equivalence class  of 
 jointly monic spans
 \begin{equation}\label{onex}\begin{diagram}A&\lTo^{c}&M&\rTo^{d}&X
 \end{diagram}.\end{equation}
 \end{defn}

 (One also says that such a $U$ a {\em family} of subobjects of $A$ {\em 
parametrized by} $X$.)  There is a partial order $\subseteq_{X}$ on 
the set of \gd subobjects of $A$ at stage $X$, defined by 
representing diagrams like (\ref{mux}).

To describe the functorality (change-of-stage) for $U\subseteq_{X}A$: given a \gd subobject $U$, represented by $(c,d)$, of $A$ at stage 
 $X$, and given $\alpha: Y\to X$, consider a diagram where the square is a 
 pull-back and the triangle is commutative
 
 \begin{equation}\label{tox}\begin{diagram}&&\cdot \SEpbk & \rTo^{d'} 
 &Y \\
&\ldTo^{c'}&\dTo&&\dTo_{\alpha }\\
A&\lTo_{c}&M& \rTo_{d}&X.
\end{diagram}\end{equation}
 Since $(c,d)$ is 
jointly monic, then so is   $(c',d')$, by an easy argument. The 
passage from $(c,d)$ to $(c',d')$ preserves the equivalence relation 
considered on the set of jointly monic spans with given ends. So we 
obtain a \gd subobject $\alpha^{*}(U)$, represented by $(c',d')$, of $A$ at stage $Y$.

Let $U\subseteq _{X}A$. If 
$\begin{diagram}Z&\rTo^{\beta}&Y&\rTo^{\alpha}&X\end{diagram}$ is a composable 
pair, then $\beta^{*}(\alpha^{*}(U))= 
(\alpha \circ \beta)^{*}(U)$ exactly, (because of the equivalence 
relation defining the notion of \gd subject; recall that a subobject 
of $A$ is not an object of $\E$).

\medskip 

There is a similar construction, but performed in the $A$-end rather than in 
the $X$-end: Given given $f:A'\to A$, and given $U\subseteq_{X}A$, 
then pulling back $c$ along $f$ provides 
a well defined subobject of $A'$ at stage $X$, which is sensibly denoted 
$f^{-1}(U)\subseteq _{X}A'$, the {\em counter image} of $U$ along $f$.
The processes $\alpha^{*}$ and $f^{-1}$ commute.

\begin{prop}Let $f:A' \to A$, and let $U'\subseteq_{X} A'$ and 
$U\subseteq _{X}A$ be \gd subobjects represented, respectively, by the 
upper and the lower spans in
\begin{equation}\label{fiix}\begin{diagram}A'&\lTo& M'&&\\
\dTo^{f}&&\dDashto_{\mu}&\rdTo
&&\\
A&\lTo&M&\rTo&X
\end{diagram}\end{equation}
Then $U'\subseteq_{X} f^{-1}(U)$ iff there exists a (necessarily unique) $\mu$ 
making the square and the triangle commute. 
\end{prop}

If further $f':A'' \to 
A'$, and $U''\subseteq _{X}A''$,  it is trivial to conclude 
a ``transitivity law'' (continuing the notation from 
the Proposition):
If $U'\subseteq_{X} f^{-1}(U)$ and $U'' \subseteq_{X} f'^{-1}(U')$, then $U'' 
\subseteq_{X} f'^{-1}(f^{-1}(U)) = (f\circ f')^{-1}(U)$.

\subsection{Partial maps  and partial 
sections $A\dashrightarrow_{U} E$}\label{ss13}
 A (generalized) {\em partial map} at stage $X$ from  
$A$ to $E$ with {\em support} $U\subseteq _{X}A$ is an equivalence class of diagrams
of the form (ignoring the $p$)
\begin{equation}\label{pmx}\begin{diagram}E&&&\\
\dTo^{p}&\luTo^{t}&&\\A&\lTo^{c}&M&\rTo^{d}&X\\
\end{diagram}\end{equation}
with $(c,d)$ jointly monic and $U$ being the \gd subobject of $A$ represented 
by $(c,d)$. The diagram represents a \gd partial {\em section} of $p$ if 
$p\circ t =c$.

The equivalence relation mentioned is the evident one, 
extending the one appearing after Definition 
\ref{AAAx}; similarly for the contravariant functorality in 
$X$. In particular, if $s$ is a \gd partial map at stage $X$, as represented by 
(\ref{pmx}), and if $\alpha:Y \to X$ is a map, then we get a \gd 
partial map $\alpha^{*}(s)$ at stage $Y$.

If $q: E\to F$ is a given actual map, one may evidently 
post-compose a \gd partial 
map $s:A\dashrightarrow _{U} E$ (with support $U\subseteq_{X}A$) with 
$q$ 
to obtain a \gd partial map $q\circ s :A \dashrightarrow _{U} F$, with 
the same support.

If $f:A'\to A$ is an actual map, we may precompose a \gd 
partial map $A\dashrightarrow _{U} E$ with $f$ to obtain a \gd partial 
map from $A'$ to $E$, with support $f^{-1}(U)\subseteq _{X} A'$.  

It is straightforward to see that the pre- and post-composition 
operations on  given \gd partial maps commute. In a full-fledged theory 
of \gd partial maps in $\E$, this will be a special case of the 
associative law for the category of \gd partial maps. We do not intend to present such a full-fledged 
theory here. It seems to require that $\E$ besides pull-backs has 
some further finite limits.

\medskip Consider spans $(c',d')$ and $(c,d)$ 
representing \gd subobjects  $U'\subseteq _{X} A'$ and $U\subseteq_{X} A$, 
respectively, both at 
stage $X$, as displayed in the right hand part of the diagram
$$
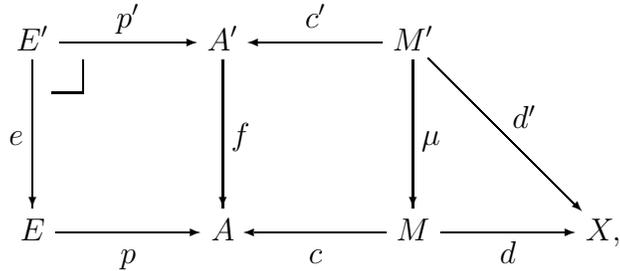
\begin{diagram}E'\SEpbk&\rTo ^{p'}& A'&\lTo^{c'}&M'&&\\
\dTo^{e} && \dTo_{f} && \dTo_{\mu}&\rdTo^{d'} &\\
E &\rTo_{p}&A&\lTo_{c}&M&\rTo_{d}&X,
\end{diagram}$$
and suppose that 
further a pull-back square is given, as in the left hand part of the 
diagram. Assume that 
$U' \subseteq _{X} f^{-1}(U)$, witnessed by the displayed map $\mu$. 
Let $t:M \to E$ represent a \gd partial section of $p$ with support 
$U$,  i.e.\  we assume that $p\circ 
t = c$.  Because of the commutativities assumed, the universal 
property of the pull-back square implies that there exists a unique 
$t': M'\to E'$ with $ p'\circ t' = c'$ and $e\circ t'=t \circ \mu$, which then represents a \gd 
partial section of $p'$ with support $U'$. Summarizing, (omitting the 
phrase ``generalized'' in a couple of places)
\begin{prop}\label{onethreex}Let $f:A'\to A$, and let $U'\subseteq _{X}A'$ and 
$U\subseteq _{X}A$ have $U'\subseteq _{X} f^{-1}(U)$. Then a partial 
section $t$ of  $p:E\to A$ with support $U$ restricts canonically to a 
partial section $t'$  of $p': E'\to A'$ with support $U'$, where $p'$ is obtained by pulling 
$p$ back along $f$.
\end{prop}
Namely, $t':M' \to E'$ is characterized by  
$p'\circ t'= c'$ and $e\circ t' = t\circ \mu$.

For $A'' \to A'\to A$, there is a rather obvious, and easily 
provable, strict associativity 
assertion which supplements this Proposition.

\subsection{Elements of a \gd subobject}\label{Elemx}
If  $(c,d)$ is a span representing a \gd subobject $U$ of $A$ at 
stage $X$, and 
$a$ is a \gd  element  of $A$ 
defined at a later stage $\alpha :Y \to X$, as depicted with full arrows in
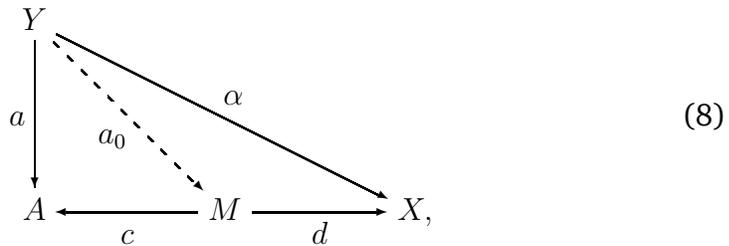
\begin{equation}\label{eltxx}\begin{diagram}Y&&&&\\
\dTo^{a}&\rdDashto_{a_{0}}\rdTo(4,2)^{\alpha}&&&\\
A&\lTo_{c}&M&\rTo_{d} &X,
\end{diagram}\end{equation}  then we say that $a\in _{\alpha} U$ if there is a 
map $a_{0}: Y \to M$ such that $c\circ a_{0}=a$ and 
$d\circ a_{0}= \alpha $. Such a map is  unique because $(c,d)$ is 
jointly monic. In this case, we say that $a_{0}$ is the {\em witness} 
or the {\em proof} (relative to the given representative $(c,d)$ for the \gd subobject 
$U$) that $a$ belongs to $U$ at the later stage.  
The assertion that $a\in _{\alpha}U$ is clearly independent of the choice 
of the span $A\leftarrow M \to X$ representing $U$.

An important special case is when no change of stage takes place, 
i.e.\ when  $Y=X$ and $\alpha :Y \to X$ is 
$id_{X}$. In this case, we say
$a\in_{X}U$. But note: $a\in_{X}U$ is not a special case of 
$a\in_{X}A$: $U$ is, unlike $A$,  not an object of $\E$, and $a\in_{X}U$ is only 
meaningful for an $a$ which is already $\in_{X}A$.

The somewhat heavy  $\in_{\alpha}$-notation can be 
dispensed with, because it is easily proved that for $U$, $\alpha$ 
and $a$ as above, 
\begin{equation}\label{mx}a\in_{\alpha} U\mbox{ iff } a\in_{Y}\alpha^{*}(U)
\end{equation}
A particular case of such $\alpha:Y\to X$ is $d:M\to X$.
We have $c \in _{d}U$, as witnessed by $id_{M}$, equivalently, we have $c\in 
_{M}d^{*}(U)$. Verbally: if $U$ is represented by $(c,d)$, then $c\in 
_{d}U$.
\medskip

It is easy to see that the relation $\in$ thus defined is stable 
under change of stage: if $a\in_{\alpha}U$, then 
$\beta^{*}(a)\in_{\alpha \circ \beta}\beta^{*}(U)$.

\subsection{Value of a partial map on an element in its 
support}\label{valx}
Let $U\subseteq _{X}A$ be represented by the span $(c,d)$ as in 
(\ref{onex}).
Consider a \gd partial map $s:A\dashrightarrow _{U}E$ with support $U$, 
thus it is  is represented by  
data as in in (\ref{pmx}) (ignoring the $p$). Consider also a \gd 
element, at the same stage $X$,  of $A$, so  
$a:X\to A$. Then if $a\in_{X} U$, witnessed by $a_{0}$ as in 
(\ref{eltxx}), then we write $s(a)$ for  $s\circ a_{0} \in _{X}E$.
More generally, if $\alpha :Y \to X$ and $a\in _{\alpha} U$, we shall 
also write $s(a)$ for the element $\in_{Y}E$ whose full notation is 
$\alpha^{*}(s)(\alpha^{*}(a))$. Experience shows that the ``change of 
of stage'' symbols $\alpha$ or $\alpha^{*}$ often can be omitted from 
notation, improving readability.

In the category of sets, this is the fundamental relation between 
{\em composition} $\circ$ in $\E$ on the one side, 
and {\em evaluation} of a function $s$  on an element $a$ on the other: 
The process leading from ``evaluation'' to ``composition $\circ$''
is the one through which we learned to compose functions - a relation 
which we want to exploit in a more general category $\E$ with 
pull-backs.

\section{Principles and constructions}

\subsection{Extensionality principle for subobjects}\label{ExtPrx}

\begin{prop}Let  $U\subseteq _{X}A$  and 
$U'\subseteq _{X}A$ be \gd subobjects of $A$. Then $U\subseteq_{X} U'$ iff 
for every $\alpha:Y\to X$ and every $a\in_{Y}A$, $a\in_{Y} 
\alpha^{*}(U)$ implies $a\in _{Y}\alpha^{*}(U')$.
\end{prop}
{\bf Proof.} Let the \gd subobject $U$ in question be represented by $(c,d)$ 
as in (\ref{onex}), and similarly let $U'$ be represented by $(c',d'): 
A \leftarrow M'\to X$. Assume that $U\subseteq _{X} U'$, witnessed by $\mu : 
M\to M'$, with the relevant commutativities. Let $\alpha:Y\to X$. 
Then if $a\in_{Y}A$ satisfies $a\in _{Y}\alpha^{*}(U)$ witnessed by 
$a_{0}:Y\to M$, it satisfies $a\in_{Y}\alpha^{*}(U')$, witnessed by $\mu\circ 
a_{0}$. 
For the converse conclusion, recall (\ref{mx}).   We may take $\alpha :Y \to X$ to be $d: 
M \to X$ and $a\in_{M}A$ to be $c$. Then 
$a\in _{M} U$ (witnessed by $id_{M}$). By assumption, therefore, $a\in _{M}U'$, and this is 
witnessed by a map $\mu :M \to M'$, which then proves the desired 
inequality between the \gd subobjects $U$ and $U'$. 

\medskip

This is a typical Kripke semantics argument, for interpreting the 
universal quantifier ``for all $a$ \ldots'', using ``for all later 
stages $Y$ and for all elements at this later stage \ldots''. Note that the ``Extensionality principle for 
subobjects'' in Proposition II.3.1 in \cite{[SDG]} only deals with actual 
subobjects, not \gd subobjects, as here.

\subsection{Yoneda principle for partial maps}
This  says roughly that 
a \gd partial map can be constructed and recognized by what it does 
to \gd elements of its 
support. 
The clue is that  both construction and recognition  are supposed to be available for 
{\em all} stages, and are preserved by change of stage. 
\begin{prop}\label{yoprx}Let $U\subseteq _{X}A$ be a \gd subobject of $A$ at stage 
$X$, represented by the span $(c,d)$, as in (\ref{onex}). Given a law 
$\sigma$ 
which to each $\alpha :Y \to X$ and  $a\in_{Y}\alpha^{*}(U)$ 
associates  $\sigma(a, \alpha)\in_{Y}E$, in a way which is stable under 
change of stage, i.e.\  $\sigma (a, \alpha)\circ \beta = \sigma(a\circ 
\beta ,\alpha \circ 
\beta)$ for any $\beta : Z\to Y$. Then there is a unique \gd partial map 
$s: A\dashrightarrow_{U}E$ such that for  
$a\in_{Y}\alpha^{*}(U)$, we have $\alpha^{*}(s)(a) = 
\sigma(a, \alpha)$.
\end{prop}
{\bf Proof}.  Take first $\alpha : Y\to X$ to be $d:M \to X$ and take 
 $a\in_{M}\alpha^{*}(U)$ to be $c$. Recall from Subsection 
\ref{Elemx} that
$c\in_{M}d^{*}(U)$. So $\sigma (c,d) \in_{M} E$. Then $s:=\sigma(c,d)$ 
is a map $M\to E$, which represents the desired \gd partial map 
$A\dashrightarrow _{U}E$. For, consider an arbitrary $\alpha:Y \to X$, and 
let $a\in_{Y} \alpha^{*}(U)$, witnessed by $a_{0}: Y\to M$. Then we 
have 
$$ s(a) = s\circ a_{0}= \sigma(c,d)\circ a_{0}= \sigma(c\circ a_{0}, d\circ a_{0}) =\sigma(a, 
\alpha )$$
using stability under change of 
state. Uniqueness of $s$ is clear.

For completeness, here is a diagram for some of the proof:
$$
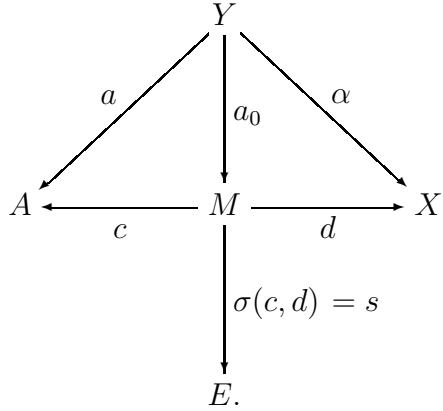
\begin{diagram}&&Y&&\\
&\ldTo^{a}&\dTo_{a_{0}}&\rdTo^{\alpha}&\\
A&\lTo_{c}&M&\rTo_{d}&X\\
&&\dTo_{\sigma(c,d)}=s&&\\
&&E.&&
\end{diagram}$$

\section{Jets associated with a relation}\label{Jex}
\subsection{Relations and their ``monads''}

 We consider the construction in (\ref{tox}), just with a change of 
 notation
 
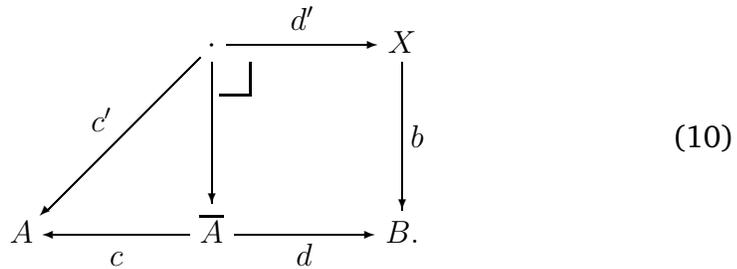
\begin{equation}\label{toxx}\begin{diagram}&&\cdot  \SEpbk & \rTo^{d'} 
 &X \\
&\ldTo^{c'}&\dTo&&\dTo_{b}\\
A&\lTo_{c}&\overline{A}& \rTo_{d}&B.
\end{diagram}\end{equation}
Then the span $(c',d')$ is jointly monic, since $(c,d)$ was assumed 
so, and hence represents a \gd subobject of $A$ at stage $X$, which 
we denote  $\M(b)$; it does not depend on the 
specific choice of span or pull-back in the construction. 
It is the ``{\em monad} around $b$'' in the 
terminology and notation from \cite{[SGM]} 2.1 and \cite{[SDG]} I.6  (where it in 
turn is borrowed from Leibniz, and from there, by non-standard analysis).  In the 
category of sets, if $b\in B$, the subset $\M(b)\subseteq A$ is 
$\{a\in A \mid (a,b)\in \M \}$.

Consider a map $\alpha :Y \to X$. It is an immediate consequence of 
the fact of ``pulling back $b\circ \alpha$ in 
stages'' that
\begin{equation}\label{Mstax}\alpha^{*}(\M(b)) = 
\M(\alpha^{*}(b)) ( = \M(b\circ \alpha)).\end{equation}
In particular, for composable $\alpha$  and $\beta$,
$$\beta^{*}(\alpha^{*}(\M(b))=(\alpha \circ \beta)^{*}(\M(b))$$
holds exactly.

\subsection{Jets and section jets}
Let $E$ be any object in $\E$, and consider a relation $\M$ from $A$ 
to $B$. 
\begin{defn} 
Let $b\in _{X}B$. 
An {\em $E$-valued jet $j$ at $b$} (relative to $\M$) is a \gd partial 
map $A\dashrightarrow _{\M(b)} 
E$.  \end{defn}
So the support of $j$ is $\M (b)\subseteq_{X}A$; the information of 
$X$ is built into $\M(b)$. Note that no specific span, like 
(\ref{relx}), is mentioned in the definition; in this sense, it is a 
``coordinate free'' definition.

Let further $p:E\to A$  be given. 
\begin{defn}\label{defsjx} Let $j$ be an 
$E$-valued jet at $b$. We say that $j$ is a {\em section jet 
of $p$ at $b$} (relative  to $\M$) if $p\circ j: A\dashrightarrow _{\M(b)}A$ is the  
partial identity map of $A$ with support $\M(b)$.
The set of such section jets, we denote $J_{\M}(b, p)$ (or $J(b,p)$ when the relation 
$\M$ is understood from the context).
\end{defn}

Note that a jet $j$, or a section jet, is not just given by 
a \gd partial map, but by 
a pair, consisting of a \gd partial map and a \gd element $b$. We shall 
mainly be concerned with section jets in the following, and it is 
worthwhile to unravel the definition in terms of representing data:

Consider 
a diagram,  where the span $(c,d)$ represents a relation $\M$ from $A$ 
to $B$, and where the right hand square is a pull-back:
$$
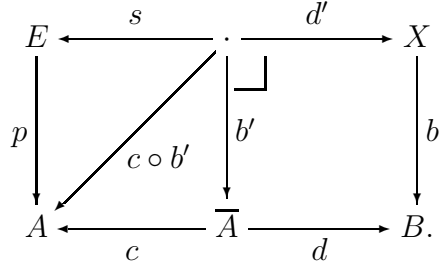
\begin{diagram}E&\lTo^{s}&\cdot \SEpbk&\rTo^{d'}& X\\
\dTo^{p}&\ldTo_{c\circ b'}&\dTo_{b'}&&\dTo_{b}\\
A&\lTo_{c}&\overline{A}&\rTo_{d}&B.
\end{diagram}$$
 The span  $(c\circ b', d'): A\leftarrow \cdot \to X$  represents the 
\gd subobject $\M(b)\subseteq _{X}A$, and the map $s$ represents an  
$E$-valued jet at $b$;  it represents a section jet if the upper 
triangle commutes, equivalently the left hand 
square, commutes.

An equivalent way of unravelling the definition diagrammatically is 
the more symmetric
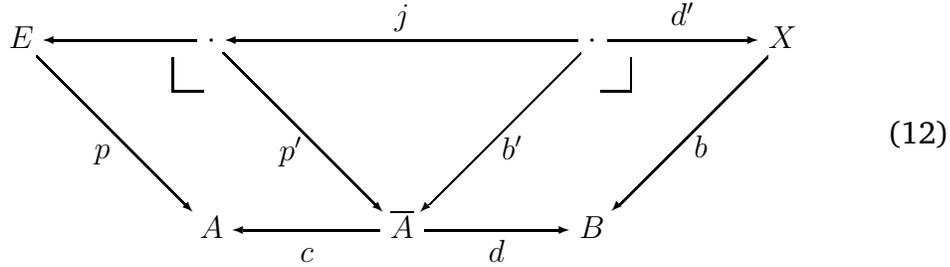
\begin{equation}\label{urx}\begin{diagram}E&\lTo & \SWpbk\cdot & 
\lTo^{j}&&&\cdot \SEpbk & \rTo^{d'} &X\\
&\rdTo_{p}&&\rdTo_{p'} & &\ldTo_{b'}  &&\ldTo_{b}\\
&&A&\lTo_{c}&\overline{A}&\rTo_{d}&B
\end{diagram}\end{equation}
where $j$ denotes the map into the (un-named) pull back, given by 
$j:=\langle s,b'\rangle$.

There is an obvious equivalence relation on the set of such 
representatives; it involves a commutativity of a square with the two 
copmarison isomorphisms between the two pull-backs involved.

So for fixed $\M$, represented by $(c,d)$, say,   and for $p:E\to A$, 
we have the set $J_{\M}(b,p)$ 
of  section jets $j$ (relative to $\M$) of $p$ at $b$. 

The set $J_{\M}(b, p)$ depends contravariantly on $b\in \E/B$ and 
covariantly   
on $p\in \E/A$. For the contravariant dependence on $b$, one sees 
from (\ref{urx}) that elements (section jets) in $J_{\M}(b, p)$ may be 
represented by factorizations $j$ of the arrow $b'$ across 
$p'$. Given a map $\alpha:Y \to X$, defining a morphism $b\circ 
\alpha \to b$ in $\E/B$, one gets the section jet $\in J_{\M}(b\circ 
\alpha )$ represented by the factorization $j\circ \alpha'$ of $b'\circ \alpha '$ 
over $p'$, where $\alpha'$ denotes some pull-back of $\alpha$ along 
$d'$. A diagrammatic rendering (which omits the pull-back square 
defining $p'$) is the following: both the displayed 
rectangles are pull-backs, and we use standard, incomplete 
notations,  e.g.\ to denote the 
object in a pull-back of $b: X\to B$  along $d$ by $d^{*}(X)$: 
\begin{equation}\label{jjjx}\begin{tikzcd}d^{*}(Y)\ar[r] \ar[d, "\alpha ' "] 
\ar[dd, bend right, shift 
right, shift right, "\alpha^{*}(j)" ']& Y\ar[d, "\alpha " ] \\
d^{*}(X) \ar[r] \ar[d, "j"]& X \ar[dd, "b"]\\
c^{*}(E) \ar[d, "p'=c^{*}(p)" ']\\
\overline{A}\ar[r, "d" ']&B.
\end{tikzcd}\end{equation}

\subsection{Morphisms between relations}

We shall now consider a morphism between two relations (between two 
different pairs of objects). Therefore, we make a change in the choice 
of letters; formerly we  considered a relation from $A$ to $B$; we 
now shall write $A_{0}$ rather than $B$. Then we have the letter $B$ 
free, and can  consider a relation $\M_{A}$ from $A$ to $A_{0}$, and another 
one $\M_{B}$ from $B$ to $B_{0}$. We are not implying that $A$ and 
$A_{0}$ are otherwise related, nor do we so for $B$ and $B_{0}$, although ultimately, 
in Section \ref{Cjx} 
on ``classical'' jet theory, where we will have $A=A_{0}$,  
$B=B_{0}$, $f=f_{0}$ (referring to (\ref{morphNx}) below).

Given a pair of maps $f:A\to B$ and $f_{0}:A_{0}\to B_{0}$, and given  relations 
$\M_{A}$ from $A$ to $A_{0}$  and $\M_{B}$ from $B$ to $B_{0}$, then there is an 
obvious notion of when $(f,f_{0})$ preserves the given relations: if the 
relations are represented by monic spans, as in 
the diagram below, this means that there is a (necessarily unique) 
map $\overline{f}$ making the two squares commute:
\begin{equation}\label{morphNx}\begin{diagram}\M_{A}:&A& \lTo& 
\overline{A}&\rTo&A_{0}\\
&\dTo^{f}&& \dTo_{\overline{f}}&&\dTo_{f_{0}}\\
\M_{B}:&B& \lTo& \overline{B}&\rTo&B_{0}.
\end{diagram}\end{equation}
In case binary products are available and chosen, this is just saying that 
$f\times f_{0}: A\times A_{0} \to B\times B_{0}$ takes the subobject 
$\M_{A}\subseteq A\imes A_{0}$ into the subobject $\M_{B}\subseteq 
B\imes B_{0}$, 
in the sense that $\M_{A} \subseteq (f\imes f_{0})^{-1}(\M_{B})$.

 Using the Extensionality Principle in Subsection \ref{ExtPrx}, it is 
 straightforward to see that this preservation property for 
 $(f,f_{0})$ can equally well be formulated in either of the 
 following ways:
 For any $a_{0}\in_{X} A_{0}$, we have
$\M_{A}(a_{0}) \subseteq f^{-1}(\M_{B}(f_{0}(a_{0}))$; or:
for any $a_{0}\in_{X}A_{0}$ and any  $a\in 
_{Y}\alpha^{*}(\M_{A}(a_{0}))$, we have $f(a)\in_{\alpha } 
\M_{B}(f_{0}(a_{0}))$.

Consider now a morphism of relations $\M_{A} \to \M_{B}$, as represented by 
the diagram (\ref{morphNx}). Let $J_{A}$ and $J_{B}$ denote the 
respective section jet constructions. Let $p:E\to B$ be given, and suppose that we have a pull-back square 
$h$ like
\begin{equation}\label{hx}\begin{diagram}
E'\SEpbk&\rTo&E\\
\dTo^{p'}&h& \dTo_{p}\\
A&\rTo_{f}&B
\end{diagram}\end{equation}
and let  $a_{0}\in _{X}A_{0}$. We shall construct a a map of sets 
\begin{equation}\label{43xx} \phi(a_{0},h): J_{B}( f_{0}(a_{0}),p ) 
\to J_{A}(a_{0}, 
p')\end{equation}
natural in $a_{0}\in \E/A_{0}$. 
An element in $J_{B}( f_{0}(a_{0}),p )$ is a section jet $j$ of $p$
$B\dashrightarrow_{\M_{B}(f_{0}(a_{0}))}E$. We  shall produce an element 
 in $J(a_{0},p')$, meaning a section jet  $j'$
$$A\dashrightarrow_{\M_{A}(a_{0})} E'.$$
We construct this \gd partial map by using the Yoneda principle in 
Proposition \ref{yoprx}: 
for any $\alpha:Y \to X$ and $a\in _{Y}A$ with $a\in_{\alpha}  \M_{A}(a_{0})$, we 
 produce an element $\sigma_{h} (a, \alpha)\in _{Y}E'$ which maps to $a$ by $p'$. 
 Because $(f,f_{0})$ 
is a morphism 
of relations from $\M_{A}$ to $\M_{B}$, we have 
 $f(a)\in_{\alpha} 
\M_{B}(f_{0}(a_{0}))$ as we observed,  so $f(a)$ belongs to the support of $j$, so $j(f(a))\in 
_{Y}E$ is 
defined as an element $\in_{Y} E$, mapping to $f(a)\in_{Y}B$ by $p$. 
By the commutativities involved, the pair $\langle a, j(f(a))\rangle $
 defines an element $\in_{Y}E'$. The law $\sigma_{h}(a, \alpha)$, given 
 by $j\mapsto \langle  a, j(f(a))\rangle $, is stable under 
change of stage $\alpha$, and so by the Yoneda principle, it 
defines the desired \gd partial map,  
$A\dashrightarrow_{\M_{A}(a_{0})} E'$.
So an explicit description of the law $\sigma_{h}$ is, for $\alpha :Y \to 
X$ and $a\in_{Y}\alpha^{*}(\M_{A}(a_{0}))$:
$$\sigma_{h}(a, \alpha ):= \langle a, j(f(a))\rangle .$$
This explicit description only mentions $j$ and $a$, so it follows 
that if we have another pull-back diagram of $f$ and $p$, say
$$\begin{diagram}
E''\SEpbk &\rTo&E\\
\dTo^{p''}&k&\dTo_{p}\\
A&\rTo_{f}&B,
\end{diagram}$$ then the resulting comparison isomorphism $\tau :E'\to E''$ 
satisfies $\tau(\sigma_{h}(a,\alpha))= \sigma_{k}(a,\alpha )$.
Therefore we not only have a map $\phi(a_{0},h)$ as in (\ref{43xx}), 
but a compatible family of maps, one for each choice of pull-back 
$f^{*}(p)$, justifying the notation
\begin{equation}\label{crx}\phi(a_{0},p): J_{B}(f_{0}(a_{0}), p) \to 
J_{A}(a_{0},f^{*}(p)).\end{equation}

Consider a composable pair of morphisms between relations $\M_{A}$, 
$\M_{B}$, and $\M_{C}$,  represented by the four right hand squares in the 
following diagram. The two further  squares $h$ and $k$ are 
assumed to be pull-backs; hence  their concatenation $k\circ h$ is likewise a 
pull-back.

\begin{equation}\label{15x}\begin{diagram}
E''\SEpbk&\rTo^{p''}&A&\lTo  & \overline{A}& \rTo &A_{0}\\
\dTo&h&\dTo_{f}&&\dTo &&\dTo_{f_{0}}\\
E'\SEpbk&\rTo^{p'}&B&\lTo &\overline{B}&\rTo&B_{0}\\
\dTo&k&\dTo_{g}&&\dTo&&\dTo_{g_{0}}\\
E&\rTo_{p}&C&\lTo&\overline{C}&\rTo& C_{0}.
\end{diagram}\end{equation}
 In these 
circumstances, we have
\begin{prop}\label{3337x}For  any $a_{0}\in_{X}A_{0}$, the composite 
$$\begin{diagram}
J_{C}(g_{0}f_{0}(a_{0}),p)&\rTo^{\phi(f_{0}(a_{0}), 
k)}&J_{B}(f_{0}(a_{0}), p')&\rTo^{\phi(a_{0},h)}&J_{A}(a_{0},p'')
\end{diagram}$$
equals $\phi(a_{0}, h\circ h')$.\end{prop}
{\bf Proof.}  Let $j: C\dashrightarrow_{U} E$ 
be an element in  the common  domain of the two maps to be compared,  with $U= 
\M_{C}(g_{0}f_{0}(a_{0}))$. The value 
of either of the two maps on this $j$ are \gd partial maps $ 
A\dashrightarrow_{\M(a_{0})} E''$.
To 
see that they are equal, we use the ``recognition'' part of the 
Yoneda priciple; is suffices to consider an arbitrary $\alpha :Y\to X$ 
and a \gd element $a\in_{Y} \alpha^{*}(\M_{A}(a_{0}))$, and see that the 
partial map 
in either case takes the value $\langle a, j(g(f(a)) \rangle$ on such 
$a$. The condition for the partial map $j$ being defined on the 
argument $g(f(a))$ follows because the relations in question 
are preserved by $(f,f_{0})$ and $(g,g_{0})$, by assumption.

\subsection{Jet bundles $\J$}\label{Jebux}

Note that if we have chosen pull-back functors $c^{*}: \E/A \to 
\E/\overline{A}$ and $d^{*}:  \E/A_{0} \to \E/\overline{A}$, (where 
$(c,d)$ is a span $A\leftarrow \overline{A}\to A_{0}$ representing $\M$) we can, using 
(\ref{urx}), 
write (for $a_{0}\in _{X}A_{0}$),
\begin{equation}\label{bijx}J_{A}(a_{0},p) \cong 
\hom_{\E/\overline{A}}(d^{*}(a_{0}),c^{*}(p))\end{equation}
(where $J_{A}$ is short for $J_{\M}$, with $\M$ is understood from the 
context).
If further $d^{*}$ admits a right adjoint $d_{*}$, we therefore also 
have
\begin{equation}\label{bij2x}{J_{A}}(a_{0},p)\cong 
\hom_{\E/A_{0}}(a_{0},d_{*}c^{*}(p)). \end{equation}

 \sloppy Consider for a fixed $p\in 
\E/A$  the contravariant set valued functor $J_{A}(-,p)$ on 
$\E/A_{0}$, and assume that it is representable (which is the case in 
(\ref{bij2x})). 
Thus there is a representing object $\J (p)={\mathbb J}_{\M}(p) \in \E/A_{0}$ 
and a bijection
$$J_{A}(a_{0},p)\cong \hom_{\E/A_{0}}(a_{0}, \J(p)).$$
More precisely, we have a 
generic section jet $\epsilon$ of $p$, defined 
at stage $\J (p)\in \E/A_{0}$, 
such that for any $a_{0} :X \to A_{0}$, and any 
$j\in J_{A}(a_{0},p)$, we 
have $\overline{j}^{*}(\epsilon)= j\in J_{A}(a_{0},p)$ for a unique 
$\overline{j}\in 
\hom_{\E/A_{0}}(a_{0}, {\mathbb J}_{\M}(p))$; $\overline{j}$ deserves the name the 
{\em classifying map} for $j$.

A rough sketch (specializing (\ref{jjjx}) of the items here (with 
incomplete but conventional  
notation, as in (\ref{jjjx})),  is found in the following 
diagram; both the displayed rectangles are a pull-backs
\begin{equation}\label{bij3x}\begin{tikzcd}d^{*}(a_{0}) \ar[r] \ar[d] \ar[dd, bend right, shift 
right, shift right, "j" ']& X\ar[d, "\overline{j}" '] 
\ar[ddd, bend left, "a_{0}"] \\
d^{*}\J(p) \ar[r] \ar[d, "\epsilon"]& \J (p) \ar[dd]\\
c^{*}(E) \ar[d, "c^{*}(p)" ']\\
\overline{A}\ar[r, "d" ']&A_{0}.
\end{tikzcd}\end{equation}
 The 
generic section jet is denoted $\epsilon$, because it, as an arrow in 
$\E/\overline{A}$, is the back 
adjunction for $d^{*}\dashv d_{*}$ (if we have the functors $d^{*}$ 
and $d_{*}$ available). And $\overline{j}$ is the classifying map for 
$j$: $\overline{j}^{*}(\epsilon)=j$, see (\ref{jjjx}), with 
different notation.

The object $\J(p)\to A_{0}$ in $\E/A_{0}$ deserves the name: the section-jet-{\em 
bundle} of $p:E\to A$ (relative to $\M$). It  exists if 
$\E$ is a locally Cartesian closed category.
The existence of $\J(p)$ may depend on $\M$, as well as on 
 $p$,  and may, if it exists, be  described, without any chosen pull-back 
functors or right adjoints for them, in terms of what \cite{[W]} calls 
{\em distributivity pull-backs}, here: a distributivity pull-back 
around 
$(c^{*}(p),d)$. The lower rectangle in (\ref{bij3x}) is a distributivity pull-back. 
In Section 5, we shall give the following verbal rendering of the 
property that makes the lower rectangle above a distributivity 
pull-back, namely: $\epsilon$ is terminal in the category of 
comorphisms over $d$ with domain $c^{*}(p)$.

If $g:A_{0}\to B_{0}$ is a map, and $q:F\to B$,  then the functor $a_{0}\mapsto 
J(g(a_{0}),q)$, 
for $a_{0}\in \E/A_{0}$ is representable as well: the pull-back 
$g^{*}{\mathbb J}(q)$ 
will do the job.  
This fact is a version of the Beck-Chevalley condition.

If every $p\in \E/A$ is provided with a (chosen)  generic  section jet 
$\epsilon$ of it, it follows 
from standard properties of bifunctors that $\J(p) \in \E/A_{0}$ depends 
functorially on $p\in \E/A$; in other words, we have (for the given 
relation $\M$ from $A$ to $B$)  a functor $$\J : 
\E/A \to \E/A_{0}:$$
{\em A relation $\M$ from $A$ to $A_{0}$ gives rise to a functor which to 
a bundle $p$ over $A$ associates its $\M$-jet-bundle $\J(p)$, which is a 
bundle over $A_{0}$.}

Given a morphism of relations $\M \to \M'$, in analogy with (\ref{43xx}), 
(denoting the corresponding jet bundle functors $\J$ and $\J '$, 
respectively). 
Then the naturality of the map described in (\ref{43xx}) implies 
that it is  mediated by a map 
\begin{equation}\label{mediatex}\Phi(h):
f^{*}(\J (p))\to  \J'(p'),\end{equation} or denoting $p'$ by 
$f^{*}(p)$, since $p'$ comes from a pull-back of $p$ along $f$,
by a map $f^{*}(\J (p))\to \J '((f^{*}(p))$
(recall the assumption that $\M$ maps into $\M'$).
In Section 
\ref{Cjx} 
below, we assume that the $\M$s are uniformly 
given, in a certain sense; we also assume  that $A=B$, $A'=B'$, and 
$f=f_{0}$, in which case
the map constructed is a map
$\Phi (f): f^{*}(\J(p))\to \J (f^{*}(p))$.

\section{Classical jets}\label{Cjx}
Our interest in relations presented by jointly monic spans 
$A\leftarrow M  \to A_{0}$ comes 
from algebraic geometry, where $A=A_{0}$ is an affine scheme, and $M$ 
is the 
``$r$th neighborhood of the diagonal'', $A_{(r)}\subseteq A\times A$ 
(or ``prolongation space $A_{(r)}$'' in \cite{[KS]} for the $C^{\infty}$ case); 
it is a reflexive symmetric relation. The role of this for jet-theory 
was made explicit in \cite{[KS]}, Chapter 1, see also \cite{[P]} IV.2; 
some of it is expounded synthetically in \cite{[SGM]}. All morphisms of affine schemes preserve the 
$r$th neighbourhood relation. 

So we consider in the following that every object $A\in \E$ 
comes together with an (endo-) relation $\M = \M_{A}$, represented by a (jointly monic) span 
\begin{equation}\label{2xx}\begin{diagram}
A&\lTo^{c}&\overline{A}&\rTo^{d}&A,
\end{diagram}\end{equation}
preserved by all maps in $\E$; in particular, for $f:A\to B$, we 
have the diagram like (\ref{morphNx}), but now with $A=A_{0}$,  $B=B_{0}$, and 
with $f=f_{0}$:
\begin{equation}\label{morphAx}\begin{diagram}A& \lTo^{c}& 
\overline{A}&\rTo^{d}&A\\
\dTo^{f}&& \dTo_{\overline{f}}&&\dTo_{f}\\
B& \lTo_{c'}& \overline{B}&\rTo_{d'}&B;
\end{diagram}\end{equation}
the map $\overline{f}$ is the witness that $f$ preserves the 
endo-relations involved. 

The jet functors for the given endo-relation $\M_{A}$ on an object 
$A$ are 
denoted $J_{A}$, and (if they exist) similarly for the jet bundles $\J_{A}$.

Because of the assumption 
 that any map 
 $f: A \to B$ preserves the given endo-relations (``takes $\overline{A}$ into 
 $\overline{B}$''),   we have in particular the 
 $\phi$-construction in (\ref{43xx}) available: 
for any   map $f:A\to B$ and for any 
 pull-back square $h$: 
 $$\begin{diagram}E'\SEpbk  &\rTo^{f_{1}}& E\\
 \dTo^{p'}&h&\dTo_{p}\\
 A¨&\rTo_{f}&B
 \end{diagram}$$
 we have, for $a_{0}\in _{X}A$, a set mapping $$\phi (a_{0},f): 
 J_{B}(f(a_{0}),p) \to 
 J_{A}(a_{0},p').$$ 
Consider two  commutative squares on top of each other
\begin{equation}\label{pb2x}
\begin{tikzcd}
F'\ar[r, "f_{2}"]\ar[d, "r' "] \ar[dd, bend right, "q' " ']
&F\ar[d, "r" ']\ar[dd, bend left, "q" ]\\
E' \ar[rd, phantom, "h"] \ar[r, "f_{1}" ]\ar[d, "p '" ]&E\ar[d, "p" ']\\
A\ar[r,"f" ']&B
\end{tikzcd}
\end{equation}
and assume that the lower square $h$ is a pull-back, and also that
the total square, which we denote $k$, is  a pull-back,  
$$\begin{tikzcd}F'\ar[rd, phantom, "k"] \ar[r, "f_{2}"] \ar[d, "q' " ']&F\ar[d, 
"q"]\\
A\ar[r, "f" ']&B
\end{tikzcd}$$
where $q'= p'\circ r'$ and $q=p\circ r$.  (Then the upper square is 
also a pull-back, but we do not need to name it). Recall by (\ref{43xx}) we 
have, for $a_{0}\in_{X}A$ a set map
$\phi(a_{0},h): J_{B}(f(a_{0}),p) \to J_{A}(a_{0},p')$, and similarly we have  
$\phi(a_{0},k): J_{B}(f(a_{0}),q) \to J_{A}(a_{0},q') $.

\begin{prop}\label{Cluex} The following diagram of sets commutes
$$\begin{diagram}J_{B}(f(a_{0}),q)&\rTo^{\phi(a_{0},k)}& J_{A}(a_{0},q')\\
\dTo^{J_{B}(f(a_{0}),r)}&&\dTo_{J_{A}(a_{0},r')}\\
J_{B}(f(a_{0}),p)&\rTo_{\phi(a_{0},h)}& J_{A}(a_{0},p').$$
\end{diagram}$$\end{prop}
{\bf Proof.} Since it is a diagram  of sets, it suffices to see that 
the value of the two composites on an element $j\in J_{B}(f(a_{0}),q)$ give 
the same element in $J_{A}(a_{0},p')$. Here, $j$ is a section jet 
$j: B\dashrightarrow_{\M (f(a_{0}))}F$ of $q$. So let 
$a\in_{\alpha}\M(a_{0})$.
Then both composites produce the element $\langle a, 
r(j(a))  \rangle \in_{Y} E'$ where $\alpha:Y \to X$. This hinges on 
a general fact, valid for generalized elements in a concatenation of 
two pull-back squares: referring to notation as in 
(\ref{pb2x}), 
it is 
the fact that 
$r'(\langle  a , c \rangle )=  \langle a  ,r(c) \rangle$ 
whenever
$a\in_{Y}A$, $c\in_{Y}F$ satisfy $f(a)= q(c)$.

\subsection{Reflexivity and symmetry}
For a relation $\M$ from $A$ to itself, it makes sense to ask whether it 
is reflexive; this is simply that it contains the diagonal relation 
$\Delta$. It is esily seen to be equivalent to: for every generalized 
element $a_{0}\in _{X}A$, we have $a_{0}\in_{X} \M(a_{0})$, - justifying the 
terminology that $\M(a_{0})$ is the monad {\em around} $a_{0}$.

Consider a reflexive relation $\M$, and a section jet $j$ at 
$a_{0}\in_{X} A$ of $p:E\to A$, so the support of $j$ is $\M(a_{0})$. Then 
since $a_{0}\in _{X}\M(a_{0})$, by reflexivity, the 
value $j(a_{0})\in _{X}E$ makes sense. 
For a section jet $j: A\dashrightarrow _{\M(a_{0})}E$ of  $p:E\to A$, we 
therefore have (cf.\ Subsection \ref{valx}) a 
particular \gd element in $E$, namely $j(a_{0})\in _{X}E$, with 
$p(j(a_{0}))=a_{0}$.

A possible symmetry of a relation $\M$ from $A$ to itself can also  
be expressed in terms of \gd elements: $a\in_{X}\M(b)$ iff 
$b\in_{X}\M(a)$.

For $\M$ reflexive and symmetric, we may think of $\M$ as encoding a 
neighbour relation: $a\M a_{0}$ meaning ``$a$ is $r$th order neighbour 
of $a_{0}$''; and $\M(a_{0})$ would then mean the $r$th order (infinitesimal) neighbourhood of 
$a_{0}$. This is how it is used in \cite{[FMSTJB]} and in \cite{[SGM]}.

 \section{Formulation in terms of fibered categories $\E^{2}$ and 
$(\E^{2})^{*}$}
Recall that if $\E$ is a category with pull-backs, one has the 
codomain fibration $\E^{2}\to \E$, associating to an object $a:X \to 
A$ its codomain $A$, and  the arrows in $\E^{2}$ are 
the commutative squares in $\E$. The fibre over $A\in \E$ thus is the 
category $\E/A$. Note that $\E$ itself is not of the form 
$\E/A$, unless we have a terminal object in $\E$.

The Cartesian arrows in $\E^{2}$ are the pull-back squares in $\E$.
Chosing pull-backs 
amounts to a clevage of the codomain fibration. Such cleavage amounts 
to giving pull-back {\em functors} $f^{*}: \E/A \to \E/A'$ for any 
$f:A'\to A$. But the present Section is ``cleavage free''.

Associated to the codomain fibration $\E^{2}\to \E$, we have its 
fibrewise dual fibration
$(\E^{2})^{*} \to \E$, as we  shall recall  below.

Some of the notions and constructions may as well be formulated for 
an arbitrary fibered category $\X\to \E$. 
This in particular 
applies to the description of the ``fibrewise dual'' fibration 
$\X^{*}\to \E$: 
\subsection{The fibrewise dual of a fibration $\X \to \E$}
Let $\pi :\X \to \E$ be any fibration. We have in mind the codomain 
fibration $\E^{2}\to \E$. In this case, $\X_{A}= \E/A$, and 
therefore, we find it convenient to denote the objects in $\X_{A}$ by 
lower case 
letters like $p$, since $p:E\to A$ is the name of a typical object in 
$\E/A$. Recall from \cite{[DFET]} or \cite{[BFF]} that an arrow in $\X^{*}$ 
over the arrow $f: A' \to A$ in $\E$, from $p' \in \X_{A'}$ to $p\in 
\X_{A}$, 
is represented by a ``vh-span'' $(v,h)$ where $v$ is vertical over $A'$ 
and $h$ Cartesian over $f$
$$\begin{diagram}\cdot& \rTo^{h}&p\\
\dTo^{v}&&\\
p'&&
\end{diagram}$$
Two such spans, say $(v,h)$ and $(v',h')$  are {\em equivalent}  if there 
is a (necessarily unique) vertical isomorphism  $i$ with $ v\circ i = v'$ 
and $ h\circ i = h'$. The equivalence classes, we call {\em 
comorphisms} from $p'$ to $p$; they are the arrows of $\X^{*}$.
The composition in $\X^{*}$ is a standard composition of spans. The functor $\X^{*}\to 
\E$ associates to  the comorphism represented by a vh span $(v,h)$, 
(as above) the arrow $\pi (h)$ in $\E$. This functor $\X^{*}\to \E$ is a fibration. A 
Cartesian arrow in $\X^{*}$ has a unique representative 
of the form $(v,h)$ with $v$ an identity arrow in $\X$.
So there is an isomorphism between the category $C(\X)$ of Cartesian arrows in 
$\X$ and the category $C(\X^{*})$ of Cartesian arrows in $\X^{*}$.
Also, $(\X^{*})_{A}$ is may be identified with $(\X_{A})^{op}$; 

Assume that we for each $A\in \E$ have an endofunctor $\J_{A}: 
\X_{A}\to \X_{A}$, hence also $\J_{A}^{op}: \X_{A}^{op}\to 
\X_{A}^{op}$. Assume also that we  have a functor $\Phi: C(\X) \to 
\X^{*}$, which for any $A \in \E$ agrees with $\J_{A}$ on  objects, 
and such that $\Phi$ agrees    with $\J_{A}^{op}$ 
on the arrows that are simultneously Cartesian and 
vertical.
Consider a vh square in $\X$, i.e.\ a commutative square in $\X$ with 
the $h_{i}$s Cartesian over $f:A'\to A$, and $w$ and $v$ vertical,
$$\begin{diagram}q_{2}&\rTo^{h_{2}}&p_{2}\\
\dTo^{w}&&\dTo_{v}\\
q_{1}&\rTo_{h_{1}}&p_{1}
\end{diagram}$$
over $f:A'\to A$ in $\E$.

There is a compatibility condition between the $\J$s and the $\Phi$,
namely commutativity of the following square in $\X^{*}$, 
$$\begin{diagram}\J_{A'}(q_{2})&\rTo^{\Phi (h_{2})}&\J_{A}(p_{2})\\
\uTo^{\J_{A'}^{op}(w)}&&\uTo_{\J_{A}^{op}(v})\\
\J_{A'}(q_{1})&\rTo_{\Phi (h)}&\J_{A}(p_{1}).
\end{diagram}$$

\begin{prop}\label{51x}Assume the compatibility for vh squares. 
Then there is a canonical functor $\J: \X^{*}\to \X^{*}$ over 
$\E$ 
agreeing with $\Phi$ on Cartesian arrows, and (except for variance) with the $\J_{A}$s on 
vertical arrows.
\end{prop}

This is an easy consequence of the way composition of comorphisms are 
defined (``standard span composition''). (The special case for $\X = 
\E^{2}$ was given  in \cite{[BFF]}, Theorem 11.1, in synthetic terms.)

\medskip The following is an explanation in fibrational terms of the notion of 
distributivity pull-back:

Given a map $d:A\to B$ in $\E$, and given $q\in \X_{A}$, we have the 
following category: its objects are comorphisms over $f$ with domain 
$q$ and codomain any object $t$ in $\X_{B}$; the arrows from the 
comorphism represented by the vh span $(e\circ v, h_{1})$ to the one represented by
$(e,h_{2})$ are represented by diagrams of the following form,
where the square is a vh-square over $d$, and where $e$ is vertical
$$\begin{diagram}\cdot & \rTo^{h_{1}} & t_{1}\\
\dTo^{v}&&\dTo\\
\cdot&\rTo_{h_{2}}&t_{2}\\
\dTo^{e}&&\\
q&&
\end{diagram}.$$
In these terms, 
if $\X\to \E$ is the codomain fibration, 
a terminal object in the category of comorphisms thus described is a distributivity pull-back in 
Weber's sense around $q,d$.

\subsection{Jets as a functor $(\E^{2})^{*}\to (\E^{2})^{*}$}

\sloppy The present Subsection gives the ultimate aim of the present note, 
namely to establish the existence of jet-bundle formation as a global 
functor $(\E^{2})^{*}\to (\E^{2})^{*}$ (for the 
``classical'' case of uniformly given endo-relations in $\E$, as in 
Section \ref{Cjx}), agreeing with the jet bundle formation $\J _{A}$ on the individual 
fibres $\E/A$ (cf.\ also \cite{[BFF]} Theorem 11.1). There are cleavage choices 
involved in the descripton of the $\J _{A}s$, and in the choice of 
right adjoints, like the $d_{*}$s used. But the text given provides a 
choice-free description of the various set-valued functors which the 
$\J_{A}$s represent, where no choice- generated coherence questions 
arise.

So consider  the case where $\X = \E^{2}$, and where each object $A$ of 
$\E$ comes with a relation $A\leftarrow \overline{A} \to A$, as in 
Section \ref{Cjx}. Therefore, we have for each $A\in \E$ a 
functor $$J_{A}:(\E/A)^{op} \times \E/A  \to Sets,$$  
with $J_{A}(b,p) $ 
as described in Definition \ref{defsjx}. We observe that 
Proposition \ref{Cluex} can be formulated in abstract fibrational terms, 
applied to  the codomain fibration 
$\E^{2} \to \E$: for, the concatenation in (\ref{pb2x}) can be seen 
as commutative square in $\E^{2}$
$$\begin{diagram}q'&\rTo^{k}&p\\
\dTo^{r'}&&\dTo_{r}\\
p'&\rTo_{h}&p
\end{diagram}$$
with the two horizontally displayed arrows being Cartesian arrows 
(namely pull-back squares in $\E$) over $f$, and the two vertically 
displayed arrows $r'$ and $r$ being vertical over $A$ and $B$, 
respectively (the equality $p\circ r =q$ is the commutativity that 
qualifies $r$ as a vertical arrow $q\to p$  over $B$, and similarly 
for $r'$). So it is a ``vh square'' over $f$. Conversely, such a vh 
square over $f$  is given by a concatenation (\ref{pb2x}) of pull-backs in 
$\E$. And in this form, the result in the Proposition can be 
interpreted as the compatibility condition needed to  establish 
existence of a functor $(\E^{2})^{*}\to  (\E^{2})^{*}$, agreeing with 
the $\J_{A}$s (except for variance) on vertical arrows, and with the 
$\Phi$ on Cartesian  arrows. Note that the value of $\Phi (h)$ as in (\ref{mediatex})
may be seen as a comorphism in $\E^{2}$ over $f$ from $\J_{A}(p') $ to 
$\J_{B}(p)$ (which in turn refers to a pull-back diagram (Cartesian 
arrow) displayed in (\ref{hx})). 
The fact that $\Phi$ preserves 
composition of Cartesian arrows follows from Proposition \ref{3337x}. 
Note that the values of $\Phi$ are not Cartesian arrows in general.

 \subsection{Jet functors in fibered categories with internal products}
 We consider a fibered category $\pi:\X \to \E$ with internal 
 products. We assume a cleavage given: for each arrow $d:M\to B$ in 
 $\E$,  we therefore have a functor $d^{*}: \X_{B}\to \X_{M}$; and we furthermore assume 
 that these functors admit right adjoints  (internal products), which 
 likewise are assumed chosen, and denoted $d_{*}$. (At this point, we are not assuming 
 Beck-Chevalley conditions.)  
 Consider a span (not assumed jointly monic) 
\begin{equation*}\begin{diagram}A&\lTo^{c}&M&\rTo^{d}&B
 \end{diagram}.\end{equation*}
 It gives rise to a functor $\J: \X_{A}\to \X_{B}$, namely the composite
 $$\begin{diagram}\X_{A}&\rTo^{c^{*}}&\X_{M}&\rTo^{d_{*}}&X_{B}.
 \end{diagram}$$
 We may also consider a similar span $(c',d')$, as displayed as the 
 upper line in the following commutative diagram
 \begin{equation}\label{morxx}\begin{diagram}A'&\lTo^{c'}&M'&\rTo^{d'}&B'\\
 \dTo^{f}&&\dTo^{\overline{f}}&h&\dTo_{g} \\
 A&\lTo^{c}&M&\rTo^{d}&B
 \end{diagram}\end{equation}
 and we similarly let $\J'$ denote $d'_{*}\circ c'^{*}: \X_{A'}\to 
 \X_{B'}$.
 Then there is a canonical natural transformation
 $g^{*}\circ \J \Rightarrow \J' \circ f^{*}: \X_{A}\to \X_{B'}$; it 
 is the 2-cell obtained by the pasting 
 $$\begin{diagram}\X_{A'}&\rTo^{c'^{*}}&\X_{M'}&\rTo^{d'_{*}}& \X_{B'}\\
 \uTo^{f^{*}}&\cong  &\uTo^{\overline{f}^{*}}&\Uparrow&\uTo_{g^{*}}\\
 \X_{A}&\rTo_{c^{*}}&\X_{M}&\rTo_{d_{*}}& \X_{B}
 \end{diagram}$$
 where the left hand 2-cell  comes from the commutativity of the 
 left hand square in (\ref{morxx}),
 and the right hand 2-cell is obtained as the mate, under $d^{*}\dashv d_{*}$, $d'^{*}\dashv 
 d'_{*}$, of the natural isomorphism
  $\overline{f}^{*}\circ d^{*}\Rightarrow d'^{*}\circ g^{*}$, which 
  in turn comes from the
  commutativity of the right hand square in (\ref{morxx}). (If this 
  square $h$ is a pull-back, the 2-cell will be an isomorphism, by the 
  Beck-Chevalley condition.)
  
  This gives a construction, in terms of the calculus of polynomial 
  functors, of the 2-cells $\Phi (h)$. There is also in these terms a 
  compatibility with the composite of two morphisms of spans, as in 
  (\ref{15x}); however, then even more canonical isomorphisms, and 
  therefore coherence questions, present themselves. 
  This is due to the choices involved in a cleavage, and of
  the choice of internal products. The approach to jet 
  bundles which we have 
  chosen in the present note was motivated by the desire to get rid 
  of such choices, by making constructions more ``coordinate free''. 
  
  A functor of the form $d_{*}\circ c^{*}$ is a special case of a 
  polynomial functor in a locally Cartesian closed category $\E$.  The possibility of a ``polynomial functor'' approach to jet bundles 
  was observed in \cite{[SGM]}, Remark 7.3.1. 
  
  This locally Cartesian closed category approach to jet bundles does not insist that the spans 
  considered $(c,d)$ are jointly monic; so they include the idea of 
  non-holonomous jets, as studied ``synthetically''  in \cite{[CNHJ]}. 
  I conjecture that the construction of $(\E^{2})^{*}\to (\E^{2})^{*}$ also 
  works for this non-holonomonous case, However, I have not been able 
  in this generality 
  to circumvent the necessity for using cleavages, i.e.\ the choice of
  pull-back functors  like $c^{*}$ and their right adjoints, and the 
  coherence questions arising.

\noindent Anders Kock,  Dept.\ of Mathematics, Aarhus Universitet, 
Denmark\\ 
\url{kock@math.au.dk}

\noindent \small Aarhus, April 2020.

\noindent \small Most diagrams were made with Paul Taylor's package.

\end{document}